\documentclass[11pt]{article}
\usepackage{amssymb}
\usepackage{amsmath}
\usepackage{amsthm}
\usepackage{latexsym}
\usepackage{amsfonts}
\usepackage{graphicx}
\usepackage{graphics}

\theoremstyle{plain}
\newtheorem{thm}{Theorem}[section]   
\newtheorem{prop}[thm]{Proposition}
\newtheorem{lem}[thm]{Lemma}

\newtheorem{Def}[thm]{Definition}

\theoremstyle{definition}

\newtheorem*{Proof}{Proof}

\newcommand{\oo} {{\omega}}

\newcommand{\bi} {{\beta}}
\newcommand{\ga} {{\gamma}}

\newcommand{\vPi} {{\varPi}}
\newcommand{\vPsi} {{\varPsi}}

\newcommand{\ld} {{\ldots}}
\newcommand{\sm} {{\smallsetminus}}
\newcommand{\thi} {{\theta}}
\newcommand{\de} {{\delta}}
\newcommand{\De} {{\varDelta}}

\newcommand{\la} {{\lambda}}
\newcommand{\el} {{\ell}}
\newcommand{\e} {{\varepsilon}}
\newcommand{\f} {{\varphi}}
\newcommand{\mi} {{\mu}}

\newcommand{\ssum}{\sum\limits}
\newcommand{\ct}{{\cal{T}}}
\newcommand{\cu}{{\cal{U}}}
\newcommand{\cp}{{\cal{P}}}
\newcommand{\ch}{{\cal{H}}}

\newcommand{\ca}{{\cal{A}}}
\newcommand{\cb}{{\cal{B}}}

\newcommand{\ra}{{\rightarrow}}
\newcommand{\oD}{{\overline{D}}}

\newcommand{\qb}{$\quad\blacksquare$}
\def\1{\it1\hspace*{-0.150cm}{\footnotesize{I}}}

\def\R{{\mathbb{R}}}
\def\C{{\mathbb{C}}}
\def\Q{{\mathbb{Q}}}
\def\N{{\mathbb{N}}}

\begin{document}
\title{\bf Common hypercyclic vectors for certain families of differential operators}
\author{\bf N. Tsirivas}

\footnotetext{{The research project is implemented within the framework of the Action ``Supporting Postdoctoral Researchers'' of the Operational Program ``Educational and Lifelong Learning'' (Action's Beneficiary: General Secretariat for Research and Technology), and is co-financed by the European Social Fund (ESF) and the Greek State.}}

\date{}
\maketitle
{\bf Abstract:} Let $(k_n)$ be a strictly increasing sequence of positive integers. If $\sum_{n=1}^{+\infty}\frac{1}{k_n}$ $=+\infty$ we establish the existence of an entire function $f$ such that for every $\lambda \in (0,+\infty )$ the set $\{ \lambda^{k_n}f^{(k_n)}(\lambda z): n=1,2,\ldots \} $ is dense in the space of entire functions endowed with the topology of uniform convergence on compact subsets of the complex plane. This provides the best possible strengthened version of a corresponding result due to Costakis and Sambarino \cite{CoSa}. From this, and using a non-trivial result of Weyl which concerns the uniform distribution modulo $1$ of certain sequences, we also derive an entire function $g$ such that for every $\lambda \in J$ the set $\{ \lambda^{k_n}g^{(k_n)}(\lambda z): n=1,2,\ldots \} $ is dense in the space of entire functions, where $J$ is ``almost" equal to the set of non-zero complex numbers. On the other hand, if $\sum_{n=1}^{+\infty}\frac{1}{k_n}<+\infty$ we show that the conclusions in the above results fail to hold.  
%
\noindent


\section{Introduction}\label{sec1}
\noindent

Let $X$ be a complex Frechet space. A sequence of continuous and linear operators $T_n:X\to X$, $n=1,2,\ldots $ is called hypercyclic provided there exists a vector $x\in X$ such that the set $\{ T_nx: n=1,2,\ldots \}$ is dense in $X$. Such a vector $x$ is called hypercyclic for the sequence $(T_n)$ and the set of hypercyclic vectors for $(T_n)$ is denoted by $\ch C(\{ T_n\})$. For a thorough study on this subject we refer to the books \cite{2}, \cite{10}.   
The symbol $\mathbb{C}$ stands for the set of complex numbers. 
We consider the set of entire functions
\[
\ch(\C):=\{f:\C\ra\C\mid f\;\text{is holomorphic}\}
\]
endowed with the topology $\ct_u$ of uniform convergence on compact subsets of $\C$. Let $K\subseteq\C$ be compact and $h:K\ra\C$ be continuous. We denote:
\[
\|h\|_K:=\max\{|h(z)|:\;z\in K\}.
\]
Let $C_n:=\oD(0,n)$ for $n=1,2,\ld$ where $D(0,n):=\{z\in\C\mid|z|<n\}$.

For $f,g\in\ch(\C)$ we set
\[
\rho(f,g):=\sum^{+\infty}_{n=1}\frac{1}{2^n}\frac{\|f-g\|_{C_n}}{1+\|f-g\|_{C_n}}.
\]
The function $\rho:\ch(\C)^2\ra\R^+$ is the usual metric of the space $(\ch(\C),\ct_u)$ as it is well known. Let $D:\ch(\C)\ra\ch(\C)$ be the usual differentiation operator, i.e. $D(f):=f'$ for $f\in\ch(\C)$, where $f'$  is the usual derivative of $f$. The operators $D^n$, $n=1,2,\ld$ acting on the space of entire functions are continuous and linear, where 
\begin{align*}
D^1:&=D \ \ \text{and} \\
D^{n+1}&=D^n\circ D \ \ \text{(the usual composition in $\ch(\C)$)}, \ \ n=1,2,\ld\;.
\end{align*}
Let $\la\in\C\sm\{0\}$ be fixed.

We consider the dilation function $\f_\la:\C\ra\C$, $\f_\la(z)=\la z$.

Now, for $n\!\!\in\!\!\N$, $\la\!\!\in\!\!\C\sm\{0\}$ we consider the linear and continuous operator\linebreak $T_{n,\la}:\ch(\C)\ra\ch(\C)$ that is defined by the following formula:
\[
T_{n,\la}(f(z)):=D^n(f\circ\f_\la(z))=\la^nf^{(n)}(\la z) \ \ f\in\ch(\C), z\in \C.
\]
It is well known that $(T_{n,\la})$ is hypercyclic. We consider the set of hypercyclic vectors for the sequence $(T_{n,\la} )$, that is 
\[
\ch C(\{T_{n,\la}\}):=\{f\in\ch(\C)\mid\overline{\{T_{n,\la}(f),\;\;n=1,2,\ld\}}=\ch(\C)\}.
\]
The set $\ch (\{T_{n,\la}\})$ is a dense, $G_\de$ subset of $(\ch (\C),T_u)$, see for instance \cite{2}, \cite{10}.

By this fact, the fact that the space $(\ch(\C),\ct_u)$ is a complete metric space and Baire's category theorem we conclude that if we have a sequence $(\la_m)$, $m=1,2,\ld$ of non-zero complex numbers, then the set $\bigcap\limits^{+\infty}_{m=1}\ch C(\{T_{n,\la_m}\})$ is a dense, $G_\de$ subset of $(\ch(\C),\ct_u)$.\ Now the following question arises naturally.

Let $I\subseteq\C\sm\{0\}$, be uncountable.\ Is it true that:
\[
\bigcap_{\la\in I}\ch C(\{T_{n,\la}\})\neq\emptyset\,?
\]
Costakis and Sambarino \cite{CoSa} showed that the above question has a positive reply for $I:=\C\sm\{0\}$ (for the greatest possible set I). Later on, Costakis \cite{7} examined some refinements of the above problem in the setting of translation operators. Following this refinement we can interpret the previous mentioned problem in the context of differential operators along \textit{sparse} powers as follows:

We fix some subsequence $(k_n)$ of natural numbers, and some positive number $\la>0$.\ We consider the set $\ch C(\{T_{k_n,\la}\})$ of hypercyclic vectors for the sequence $\{T_{k_n,\la}\}$.\ That is
\[
\ch C(\{T_{k_n,\la}\}):=\{f\in\ch(\C)\mid\overline{\{T_{k_n,\la}(f),\;\;n=1,2,\ld\}}=\ch(\C)\}.
\]
Let $I\subseteq\C\sm\{0\}$, be uncountable.\ Is it true that:
\[
\bigcap_{\la\in I}\ch C(\{T_{k_n,\la}\})\neq\emptyset\,?
\]
In this direction we prove in Section \ref{sec2} the following proposition.
\begin{prop}\label{prop1.1}
Let $(k_n)$, $n=1,2,\ld$ be a subsequence of natural numbers such that:
\[
\sum^{+\infty}_{n=1}\frac{1}{k_n}=+\infty.
\]
Then the set
\[
\ch:=\bigcap_{\la\in(0,+\infty)}\ch C(\{T_{k_n,\la}\})
\]
is residual in $(\ch(\C),\ct_u)$, i.e. it contains a dense $G_{\de}$ set.
\end{prop}

Further, we can ask the following question:

Is it true that
\[
\ch:=\bigcap_{z\in\C\sm\{0\}}\ch C(\{T_{k_n,z}\})\neq\emptyset\,?
\]
We do not know the answer to this question.
However, we can give a positive reply for a set $J\subseteq\C\sm\{0\}$ that is nearly equal to $\C\sm\{0\}$. For the sequel we denote by $\la_1$ the usual Lebesgue measure on the real line. In order to describe the set $J$ we just mentioned, let us give the following definitions:
\begin{Def}\label{Def1.2}
Let some subset $J\subseteq\C\sm\{0\}$ and some positive number $\la>0$.

We consider the set $J_\la\subseteq[0,1)$, defined by:
\[
J_\la:=\{\thi\in[0,1)\mid\exists\;z\in J:\;z=\la\cdot e^{2\pi i\thi}\},
\]
and we call this set the set of arcs of $J$ with respect to $\la$.

Of course, for every subset $J\subseteq\C\sm\{0\}$ the sets $J_\la$ are well defined, but we may have $J_\la=\emptyset$ for some $\la\in(0,+\infty)$.
\end{Def}
\begin{Def}\label{Def1.3}
A subset $J\subseteq\C\sm\{0\}$ is called of full measure if for every $\la>0$ we have $\la_1(J_\la)=1$.
\end{Def}

After the above definitions in Section \ref{sec3} we prove the following Theorem \ref{thm1.4}.
\begin{thm}\label{thm1.4}
Let $(k_n)$ be a subsequence of natural numbers such that $\ssum^{+\infty}_{n=1}\dfrac{1}{k_n }=+\infty$.\ Then there exists some full measure set $J\subseteq\C\sm\{0\}$ such that:
\[
\bigcap_{z\in J}\ch C(\{T_{k_n,z}\})\neq\emptyset.
\]
\end{thm}

We mention that given a sequence $(k_n)$ of positive integers such that the series $\sum_{n=1}^{+\infty}\frac{1}{k_n}$ converges, then for any interval $I$ of the positive (or negative) real line we have $\bigcap_{\lambda \in I}\ch C(\{T_{k_n,\lambda }\})=\emptyset$. We do not give a proof of this result. Instead, we refer to the proof of item $(i)$ of Theorem 1.1 in \cite{Tsi4}, which can be easily adapted to our case. The interested reader will have no trouble to check the details.

The above theorem, Theorem \ref{thm1.4}, is the main result of this paper and its proof uses Proposition \ref{prop1.1} and a result of Weyl on uniform distribution of sequences. The proof of Proposition \ref{prop1.1} refines the argument in the proof of the common hypercyclicity criterion from \cite{CoSa}. A similar approach has been recently developed by the author in \cite{Tsi4}. In \cite{Tsi4}, we deal with certain families of backward shift operators. Let us mention that new common hypercyclicity  criteria are available, due to the work of Bayart and Matheron \cite{3} and Shkarin \cite{13}. Unfortunately, it is not clear to us if these criteria are applicable in our case. Similar problems for translation type operators are treated in a series of recent papers, \cite{Ba}, \cite{8}, \cite{Tsi}, \cite{Tsi2}, \cite{Tsi3}, which require a different approach at a technical level. For a sample of results on common hypercyclic vectors, see \cite{AbGo}, \cite{3}-\cite{5}, \cite{LeMu}-\cite{13}.
\section{The solution of our problem in the basic case}\label{sec2}
\noindent

In order to prove Proposition \ref{prop1.1} we need some preparation. Let $\vPsi:=\{p_1,p_2,\ld\}$ be an enumertaion of all non-zero complex polynomials with coefficients in $\Q+i\Q$.

For every $\rho=2,3,\ld$ we consider the set $\De_\rho:=\Big[\dfrac{1}{\rho},\rho\Big]$.\ Obviously, $\bigcup\limits^{+\infty}_{\rho=2}\De_\rho=(0,+\infty)$.\ From now on we fix some subsequence $(k_n)$ of natural numbers.

For every $n$, $j$, $m\in\N$, $\rho>2$, $s>1$ we consider the set
\[
E(n,\rho,j,s,m):=\bigg\{f\in\ch(\C)\mid\forall\;\la\in\De_\rho\;\exists\;v\le m:
\|T_{k_v,\la}(f)-p_j\|_{C_n}<\frac{1}{s}\bigg\}.
\]
\begin{lem}\label{lem2.1}
For every $n$, $j$, $m\in\N$ and $\rho$, $s>2$ the set $E(n,\rho,j,s,m)$ is open in $(\ch(\C),\ct_u)$.
\end{lem}
\begin{lem}\label{lem2.2}
We set
\[
G:=\bigcap^{+\infty}_{n=1}\bigcap^{+\infty}_{\rho=1}\bigcap^{+\infty}_{j=1}
\bigcap^{+\infty}_{s=1}\bigcup^{+\infty}_{m=1}E(n,\rho,j,s,m).
\]
Then
\[
G\subseteq\ch.
\]
\end{lem}
\begin{lem}\label{lem2.3}
For every $n$, $j$, $\in\N$ and $\rho,s>2$, the set $\bigcup^{+\infty}_{m=1}E(n,\rho,j,s,m)$ is dense in $(\ch(\C),\ct_u)$.
\end{lem}

The proof of the above lemmas \ref{lem2.1} and \ref{lem2.2} are left to the interested reader We will prove later only Lemma \ref{lem2.3}.
\begin{lem}\label{lem2.4}
Let some fixed $m_0\in\N$, $\la_0>0$ and a non-zero polynomial $p(z)=\ssum^{\el_0}_{j=0}\bi_jz^j$ for $\bi_j\in\C$, $j=0,1,\ld,\el_0$, $\el_0\in\N\cup\{0\}$.\ We consider the polynomial
\[
f(z):=\sum^{\el_0}_{j=0}\frac{j!}{(j+m_0)!}\frac{\bi_j}{\la^{j+m_0}_0}z^{j+m_0}.
\]
Then the above polynomial consists a solution of the differential equation
\[
T_{m_0,\la_0}(y)=p.  \eqno{(\ast)}
\]
\end{lem}
\begin{Proof}
Let $f(z)=\ssum^{+\infty}_{n=0}c_nz^n$ be a power series with center 0 and radius of convergence $+\infty$, $c_n\in\C$ for $n=0,1,2,\ld\;.$ Then for the derivative $f^{(k)}$ of $f$ for $k=1,2,\ld$ it holds:
\[
f^{(k)}(z)=\sum^{+\infty}_{n=k}n(n-1)\times\cdots\times(n-k+1)c_nz^{n-k}=
\sum^{+\infty}_{n=0}(n+1)(n+2)\times\cdots\times(n+k)c_{n+k}z^n \leqno{\text{(i)}}
\]
for every $z\in\C$ and
\[
f^{(k)}(0)=k!c_k \ \ \text{for} \ \ k=0,1,2,\ld\;. \leqno{\text{(ii)}}
\]
We apply the previous for a polynomial.

Let us take a polynomial $f(z)\neq0$, $f(z)=\ssum^{N_0}_{v=0}a_vz^v$ for some $N_0\in\N$, $N_0\ge m_0$ $a_v\in\C$, $v=0,1,\ld,N_0$.

For $1\le m\le N_0$ we have:
\begin{align}
f^{(m)}(z)&=\sum^{N_0-m}_{n=0}(n+1)(n+2)\times\cdots\times(n+m)c_{n+m}z^n \nonumber\\
&=\sum^{N_0-m}_{n=0}(n+1)(n+2)\times\cdots\times(n+m)\frac{f^{(n+m)}(0)}{(n+m)!}z^n \nonumber\\
&=\sum^{N_0-m}_{n=0}\frac{f^{(n+m)}(0)}{n!}z^n, \ \ \text{for every} \ \ z\in\C. \label{eq1}
\end{align}
Let $m_0\in\N$, $\la_0>0$ as in the hypothesis.\ For $z\in\C$, we compute:
\begin{align}
(T_{m_0,\la_0}(f))(z)&=\la^{m_0}_0f^{(m_0)}(\la_0z)\overset{(1)}{=}
\la^{m_0}_0\cdot\sum^{N_0-m_0}_{n=0}\frac{f^{(n+m_0)}(0)}{n!}(\la_0z)^n \nonumber\\
&=\sum^{N_0-m_0}_{n=0}\frac{\la^{n+m_0}_0f^{(n+m_0)}(0)}{n!}z^n.  \label{eq2}
\end{align}
We suppose now that a polynomial $f\neq0$ with $\deg f(z)=N_0\ge m_0$ is a solution of differential equation $(\ast)$.\ Then $T_{m_0,\la_0}(f)=p$.\ By (\ref{eq2}) we get:
\begin{eqnarray}
\sum^{N_0-m_0}_{n=0}\frac{\la^{n+m_0}_0f^{(n+m_0)}(0)}{n!}z^n=
\sum^{\el_0}_{j=0}\bi_jz^j.  \label{eq3}
\end{eqnarray}
By (\ref{eq3}) we take:

1)\; $N_0=\el_0+m_0$ and \medskip

2) \; $\dfrac{\la^{n+m_0}_0f^{(n+m_0)}(0)}{n!}=\bi_n$ for every $n=0,1,\ld,\el_0$.

The previous equality 2) gives us that:
\[
\frac{f^{(n+m_0)}(0)}{(n+m_0)!}=\frac{n!\bi_n}{(n+m_0)!\cdot\la^{n+m_0}_0} \ \ \text{for every} \ \ n=0,1,\ld,\el_0.
\]
This gives that the polynomial
\[
f(z)=\sum^{\el_0}_{j=0}\frac{j!}{(j+m_0)!}\frac{\bi_j}{\la^{j+m_0}_0}z^{j+m_0}
\]
is a possible solution of differential equation $(\ast)$.

Indeed! with a straightforward computation we prove now that the above polynomial $f(z)$ is a solution of the differential equation $(\ast)$.

We have:
\begin{align*}
T_{m_0,\la_0}(f)(z)&=\la^{m_0}_0f^{(m_0)}(\la_0z)\\
&=\la^{m_0}_0\cdot\sum^{N_0-m_0}_{n=0}(n+1)(n+2)
\times\cdots\times(n+m_0)
\frac{f^{(n+m_0)}(0)}{(n+m_0)!}(\la_0z)^n\\
&=\la^{m_0}_0\cdot\sum^{N_0-m_0}_{j=0}(j+1)(j+2)\times\cdots\times(j+m_0)\frac{j!}
{(j+m_0)!}\frac{\bi_j}{\la^{j+m_0}_0}\la^j_0z^j \\
&=\sum^{\el_0}_{j=0}\bi_jz^j=p(z).
\end{align*}

So, polynomial $f(z)$ solves the differential equation $(\ast)$. \qb
\end{Proof}

Now, let fixed $m_0\in\N$, $\la_0>0$ and polynomial $p(z)=\ssum^{\el_0}_{j=0}\bi_jz^j$, $\el_0=\deg p(z)$ $p\neq0$.

The polynomial
\[
f(z)=\sum^{\el_0}_{j=0}\frac{j!}{(j+m_0)!}\frac{\bi_j}{\la^{j+m_0}_0}z^{j+m_0}
\]
(that is a solution of the differential equation: $T_{m_0,\la_0}(y)=p$) is called the solution with data $(m_0,\la_0,p)$ of the previous differential equation and we say simply the solution $(m_0,\la_0,p)$ of the equation $(T_{m_0,\la_0}(y)=p$, or the solution $(m_0,\la_0,p)$.
\begin{lem}\label{lem2.5}
Let fixed $m_0\in\N$, $\la_0>0$ and a non-zero polynomial $p(z)=\ssum^{\el_0}_{j=0}\bi_jz^j$, $\el_0\in\{0,1,2,\ld\}$.\ Let $f$ be the solution $(m_0,\la_0,p)$ of the differential equation
\[
T_{m_0,\la_0}(y)=p. \eqno(\ast)
\]
We set $N_0:=\deg f(z)=m_0+\el_0$, $\el_0=\deg p(z)$.\ We fix some positive number $\e_0\in(0,1)$ and some positive number $R_0>1$.

We set $M_0:=\max\{|\bi_j|,\;j=0,1,\ld,\el_0\}$, $M_1:=M_0\cdot\ssum^{\el_0}_{j=0}R^j_0$.

Then, for every positive number $\la\in\big[\la_0,\la_0\cdot\sqrt[N_0]{1+\e_0/M_1}\big)$ the following inequality holds:
\[
\|T_{m_0,\la_0}(f)-T_{m_0,\la}(f)\|_{\oD_{R_0}}<\e_0,
\]
where $\oD_{R_0}=\oD(0,R_0):=\{z\in\C\mid|z|\le R_0\}$.
\end{lem}
\begin{Proof}
By relation (\ref{eq2}) of Lemma \ref{lem2.4} we have:
\setcounter{equation}{0}
\begin{eqnarray}
(T_{m_0,\la_0}(f))(z)=\sum^{N_0-m_0}_{n=0}\frac{\la^{n+m_0}_0f^{(n+m_0)}(0)}
{n!}z^n.  \label{eq1}
\end{eqnarray}
By relation (\ref{eq1}) of Lemma \ref{lem2.4} we have for some $\la>\la_0$, $z\in\C$.
\begin{eqnarray}
f^{(m_0)}(\la z)=\sum^{N_0-m_0}_{n=0}\frac{f^{(n+m_0)}(0)}{n!}\la^nz^n.  \label{eq2}
\end{eqnarray}
So, we have by (\ref{eq2}), for $\la>\la_0$, $z\in\C$
\begin{align}
T_{m_0,\la}(f)(z)&=\la^{m_0}f^{(m_0)}(\la z)=\la^{m_0}\cdot\sum^{N_0-m_0}_{n=0}\frac{f^{(n+m_0)}(0)}{n!}\la^n z^n \nonumber \\
&=\sum^{N_0-m_0}_{n=0}\frac{f^{(n+m_0)}(0)}{n!}\la^{n+m_0}z^n.  \label{eq3}
\end{align}
Thus, by (\ref{eq1}) and (\ref{eq3}), for $\la>\la_0$ and $z\in\C$ we get:
\begin{align}
|(T_{m_0,\la_0}(f))(z)-(T_{m_0,\la}(f))(z)|=&\,\bigg|\sum^{N_0-m_0}_{n=0}
\la^{n+m_0}_0\frac{f^{(n+m_0)}(0)}{n!}z^n\nonumber \\ &-\sum^{N_0-m_0}_{n=0}\frac{f^{(n+m_0)}(0)}{n!}
\la^{n+m_0}z^n\bigg|\nonumber \\
=&\,\bigg|\sum^{\el_0}_{n=0}(\la_0^{n+m_0}-\la^{n+m_0})\frac{f^{(n+m_0)}(0)}{n!}z^n\bigg|\nonumber \\
\le&\,\sum^{\el_0}_{n=0}|\la^{n+m_0}-\la^{n+m_0}_0|\bigg|\frac{f^{(n+m_0)}(0)}{n!}\bigg|
|z|^n \nonumber \\
=&\,\sum^{\el_0}_{n=0}\bigg|\bigg(\frac{\la}{\la_0}\bigg)^{n+m_0}-1\bigg|\cdot
\bigg|\frac{\la^{n+m_0}_0f^{(n+m_0)}(0)}{n!}\bigg||z|^n.  \label{eq4}
\end{align}
But $\bi_n:=\dfrac{\la^{n+m_0}_0f^{(n+m_0)}(0)}{n!}$ for every $n=0,1,\ld,\el_0$ because $f$ is the solution $(m_0,\la_0,p)$.

So, by (\ref{eq4}) we get for $z\in\C$, $\la>\la_0$
\begin{align}
|(T_{m_0,\la_0}(f))(z)-(T_{m_0,\la}(f))(z)|&\le\sum^{\el_0}_{n=0}\bigg|
\bigg(\frac{\la}{\la_0}\bigg)^{n+m_0}-1\bigg||\bi_n|\,|z|^n \nonumber \\
&\le\sum^{\el_0}_{n=0}\bigg(\bigg(\frac{\la}{\la_0}\bigg)^{N_0}-1\bigg)M_0|z|^n \nonumber \\
&=\bigg(\bigg(\frac{\la}{\la_0}\bigg)^{N_0}-1\bigg)\cdot M_0\cdot\sum^{\el_0}_{n=0}|z|^n.  \label{eq5}
\end{align}
For $z\in\oD_{R_0}$ we have by (\ref{eq5})
\begin{align}
|T_{m_0,\la_0}(f)(z)-T_{m_0,\la}(f)(z)|&\le M_0\cdot\sum^{\el_0}_{n=0}R^n_0\cdot\bigg(\bigg(\frac{\la}{\la_0}\bigg)^{N_0}-1\bigg) \nonumber \\
&=M_1\bigg(\bigg(\frac{\la}{\la_0}\bigg)^{N_0}-1\bigg).  \label{eq6}
\end{align}
Inequality (\ref{eq6}) holds for every positive $\la\ge\la_0$.\ So, for $\la\in\big[\la_0,\la_0\cdot\sqrt[N_0]{1+\e_0/M_1}\big)$ we have:
\begin{align}
\la<\la_0\cdot\sqrt[N_0]{1+\dfrac{\e_0}{M_1}}&\Leftrightarrow\frac{\la}{\la_0}<
\sqrt[N_0]{1+\dfrac{\e_0}{M_1}}\Leftrightarrow\bigg(\frac{\la}{\la_0}\bigg)^{N_0}<1+
\frac{\e_0}{M_1}  \nonumber \\
&\Leftrightarrow\bigg(\frac{\la}{\la_0}\bigg)^{N_0}-1<\frac{\e_0}{M_1}\Leftrightarrow
M_1\cdot\bigg(\bigg(\frac{\la}{\la_0}\bigg)^{N_0}-1\bigg)<\e_0.  \label{eq7}
\end{align}
By (\ref{eq6}), (\ref{eq7}) and the fact that the function $T_{m_0,\la_0}(f)-T_{m_0,\la}(f)$ is a polynomial and so continuous we get finally:
\[
\|T_{m_0,\la}(f)-T_{m_0,\la_0}(f)\|_{\oD_{R_0}}<\e_0 \ \ \text{for every} \ \ \la\in\big[\la_0,\la_0\cdot\sqrt[N_0]{1+\dfrac{\e_0}{M_1}}\big)
\]
and the proof of this lemma is complete. \qb
\end{Proof}

For the following lemma we consider some data.

More specifically:

Let two fixed polynomials $p$ and $Q$ where $p(z)=\ssum^{\el_0}_{j=0}\bi_jz^j$, $p\neq0$ for $\bi_i\in\C$, $i=0,1,\ld,\el_0$, $\el_0\in\{0,1,2,\ld\}$.

Let fixed $R_0\in(1,+\infty)$.\ We set
\[
M_0:=\max\{|\bi_j|\mid\;j=0,1,\ld,\el_0\}.
\]
We consider the sequence
\[
\ga_v:=\frac{(2R_0)^v}{v!}, \ \ v=1,2,\ld\;.
\]
We prove that $\ga_v\ra0$.

Really, by the ratio criterion for convergence of sequences we get:
\[
\frac{\ga_{v+1}}{\ga_v}\equiv\frac{\dfrac{(2R_0)^{v+1}}{(v+1)!}}{\dfrac{(2R_0)^v}{v!}}=
\frac{2R_0}{v+1}\ra0 \ \ \text{as} \ \ v\ra+\infty, \ \ \text{that gives} \ \ \ga_v\ra0.
\]
So, there exists some natural number $N_0\!\!\in\!\!\N$ such that: (because $M_0\cdot\el_0!$ $\ga_v\ra0$)\linebreak $M_0\,\el_0!$ $\ga_v<1$ for every $v\in\N$, $v\ge N_0$.

We set $N_1:=\max\{N_0,\deg Q,\el_0\}+1$.

We consider now two positive fixed numbers $a_0,b_0$, where $0<a_0<1<b_0<+\infty$, some fixed natural number $v_0>2$ and some partition $\De\!=\!\{a_0\!=\!\de_1\!<\!\de_2\!<\!\cdots\!<\!\de_{v_0}\!=\!b_0\}$ of the closed interval $[a_0,b_0]$.

We consider now some fixed finite sequence of natural numbers $m_1,m_2,\ld,m_{v_0}$, such that $m_1<m_2<\cdots<m_{v_0}$ and $m_1>N_1$, and $m_{i+1}-m_i>N_1$, $i=1,2,\ld,v_0-1$.\ For every $i=1,2,\ld,v_0$ we consider the solution $(m_i,\de_i,p)$ that we denote by $f_i$, $i=1,2,\ld,v_0$ for simplicity.\ We denote
\[
\vPi:=\sum^{v_0}_{i=1}f_i+Q.
\]
By the above notations and terminology, we have the following lemma:
\begin{lem}\label{lem2.6}
Let arbitrary $\la\in[a_0,b_0]$.\ We consider the unique $i\in\{1,2,\ld,v_0-1\}$ such that $\la\in[\de_i,\de_{i+1})$, (if there exists, or else $\la=b_0$).\ Then the following inequality holds:
\[
\|T_{m_i,\la}(\vPi)-p\|_{\oD_{R_0}}<\|T_{m_i,\de_i}(f_i)-T_{m_i,\la}(f_i)\|_{\oD_{R_0}}
+\frac{1}{2^{m_{i+1}-(m_i+2)}}, \eqno{(\ast)}
\]
if $\la\in[a_0,b_0)$ or else
\[
\|T_{m_{i_0},b_0}(\vPi)-p\|_{\oD_{R_0}}=0 \ \ \text{if} \ \ \la=b_0.
\]
\end{lem}
\begin{Proof}
The case $\la_0=b_0$ is obvious.\ Let $\la_0\in[a_0,b_0)$, be fixed.\ We consider the unique $i_0\in\{1,2,\ld,v_0-1\}$ such that $\la_0\in[\de_{i_0},\de_{i_0+1})$.

We have:
\[
\vPi:=Q+\sum^{v_0}_{i=1}f_i \ \ \text{by definition}.
\]
The operator $T_{m_{i_0},\la_0}$ is linear.\ So, we compute:
\setcounter{equation}{0}
\begin{align}
T_{m_{i_0},\la_0}(\vPi)&=T_{m_{i_0},\la_0}\bigg(Q+\sum^{v_0}_{i=1}f_i\bigg) \nonumber\\
&=T_{m_{i_0},\la_0}(Q)+\sum^{v_0}_{i=1}T_{m_{i_0},\la_0}(f_i).  \label{eq1}
\end{align}
As it is well known we have:
\begin{eqnarray}
T_{m_{i_0},\la_0}(Q)=\la^{m_{i_0}}_0Q^{(m_{i_0})}(\la_0z)=0  \label{eq2}
\end{eqnarray}
because $\deg Q<N_1<m_1\le m_{i_0}$ by our hypothesis.

Similarly, if $i_0>1$ we have for $1\le i\le i_0-1$
\begin{eqnarray}
T_{m_{i_0},\la_0}(f_i)=\la^{m_{i_0}}_0f^{(m_{i_0})}_i(\la_0z)=0  \label{eq3}
\end{eqnarray}
because $\deg f_i=m_i+\el_0\le m_{i_0-1}+\el_0<m_{i_0}$ by the definition of the natural number $N_1$ in our hypothesis.

By (\ref{eq1}), (\ref{eq2}) and (\ref{eq3}) we get:
\begin{eqnarray}
T_{m_{i_0},\la_0}(\vPi)=\sum^{v_0}_{i=i_0}T_{m_{i_0},\la_0}(f_i).  \label{eq4}
\end{eqnarray}
We have now:
\begin{align}
\sum^{v_0}_{i=i_0}T_{m_{i_0},\la_0}(f_i)=&\,T_{m_{i_0},\la_0}(f_{i_0})+\sum^{v_0}_{i=i_0+1}
T_{m_{i_0},\la_0}(f_i) \nonumber\\
=&\,(T_{m_{i_0},\la_0}(f_{i_0})-T_{m_{i_0},\de_{i_0}}(f_{i_0}))+T_{m_{i_0},\de_{i_0}}(f_{i_0})\nonumber \\
&+\sum^{v_0}_{i=i_0+1}T_{m_{i_0},\la_0}(f_i).  \label{eq5}
\end{align}
By (\ref{eq4}) and (\ref{eq5}) we get now:
\begin{align}
T_{m_{i_0},\la_0}(\vPi)-p=&\,(T_{m_{i_0},\la_0}(f_{i_0})-T_{m_{i_0},\de_{i_0}}(f_{i_0}))\nonumber \\
&\,+(T_{m_{i_0},\de_{i_0}}(f_{i_0})-p)+\sum^{v_0}_{i=i_0+1}T_{m_{i_0},\la_0}(f_i)\nonumber\\ &\,=(T_{m_{i_0},\la_0}(f_{i_0})-T_{m_{i_0},\de_{i_0}}(f_{i_0}))+\sum^{v_0}_{i=i_0+1}T_{m_{i_0},\la_0}(f_i) \label{eq6}
\end{align}
because $f_{i_0}$ is the $(m_{i_0},\de_{i_0},p)$ solution.

By (\ref{eq6}) and the triangle inequality we get:
\begin{align}
\|T_{m_{i_0},\la_0}(\vPi)-p\|_{\oD_{R_0}}\le&\,\|T_{m_{i_0}},\de_{i_0}(f_{i_0})-T_{m_{i_0},\la_0}
(f_{i_0})\|_{\oD_{R_0}}\nonumber \\
&+\sum^{v_0}_{i=i_0+1}\|T_{m_{i_0},\la_0}(f_i)\|_{\oD_{R_0}}.  \label{eq7}
\end{align}
So by (\ref{eq7}) we have to estimate the quantities $\|T_{m_{i_0},\la_0}(f_i)\|_{\oD_{R_0}}$ for $i=i_0+1,\ld,v_0$.

So, we fix some $j_0\in\{i_0+1,\ld,v_0\}$. We estimate the quantity:
\[
\|T_{m_{i_0},\la_0}(f_{j_0})\|_{\oD_{R_0}}.
\]

We see easily that for the $m_{i_0}$-th $(z^{k+m_{j_0}})^{(m_{i_0})}$ derivative of the monomial $z^{k+j_0}$ holds inductively:
\begin{eqnarray}
(z^{k+m_{j_0}})^{(m_{i_0})}=(k+m_{j_0})\cdot(k+m_{j_0}-1)\times\cdots\times
(k+m_{j_0}-m_{i_0}+1)\cdot z^{k+m_{j_0}-m_{i_0}},  \label{eq8}
\end{eqnarray}
for every $k=0,1,\ld,\el_0$.

The polynomial $f_{j_0}$ is the $(m_{i_0},\de_{i_0},p)$ solution so
\[
f_{j_0}(z)=\sum^{\el_0}_{k=0}\frac{k!}{(k+m_{j_0})!}\frac{\bi_k}{\de^{k+m_{j_0}}_{j_0}}
z^{k+m_{j_0}}.
\]
Thus, by (\ref{eq8}) we get for the $m_{i_0}$-derivative of $f_{j_0}(z)$.
\begin{align}
f^{(m_{i_0})}_{j_0}(z)=&\,\bigg(\sum^{\el_0}_{k=0}\frac{k!}{(k+m_{j_0})!}
\frac{\bi_k}{\de^{k+m_{j_0}}_{j_0}}z^{k+m_{j_0}}\bigg)^{(m_{i_0})} \nonumber \\
=&\,\sum^{\el_0}_{k=0}\frac{k!}{(k+m_{j_0})!}\frac{\bi_k}{\de^{k+m_{j_0}}_{j_0}}
(z^{k+m_{j_0}})^{(m_{i_0})} \nonumber \\
=&\,\sum^{\el_0}_{k=0}\frac{k!}{(k+m_{j_0})!}\frac{\bi_k}{\de^{k+m_{j_0}}_{j_0}}(k+m_{j_0})(k+m_{j_0}-1)\nonumber \\
&\times\cdots\times(k+m_{j_0}-m_{i_0}+1)\cdot z^{k+m_{j_0}-m_{i_0}} \nonumber\\
=&\,\sum^{\el_0}_{k=0}\frac{k!\bi_k}{\de^{k+m_{j_0}}_{j_0}}\frac{1}
{(k+m_{j_0}-m_{i_0})!}\cdot z^{k+m_{j_0}-m_{i_0}}.  \label{eq9}
\end{align}
Using (\ref{eq9}) we compute: for $z\in\C$
\begin{align}
(T_{m_{i_0},\la_0}(f_{j_0}))(z)&=\la_0^{m_{i_0}}f_{j_0}^{(m_{i_0})}(\la_0z) \nonumber \\
&=\la_0^{m_{i_0}}\cdot\sum^{\el_0}_{k=0}\frac{k!\bi_k}{\de^{k+m_{j_0}}_{j_0}}\cdot\frac{1}{(k+m_{j_0}-m_{i_0})!}\cdot
\la^{k+m_{j_0}-m_{i_0}}_0z^{k+m_{j_0}-m_{i_0}}\nonumber \\
&=\sum^{\el_0}_{k=0}k!\bi_k\bigg(\frac{\la_0}{\de_{j_0}}\bigg)^{k+m_{j_0}}\cdot\frac{1}{(k+m_{j_0}-m_{i_0})!}
z^{k+m_{j_0}-m_{i_0}}.  \label{eq10}
\end{align}
We have $\la_0<\de_{j_0}$.\ So by (\ref{eq10}) for $z\in\oD_{R_0}$ we have:
\begin{align}
|T_{m_{i_0},\la_0}(f_{j_0})(z)|&\le\sum^{\el_0}_{k=0}\bigg|k!\bi_k\cdot\frac{1}
{(k+m_{j_0}-m_{i_0})!}\bigg|\cdot R^{k+m_{j_0}-m_{i_0}}_0\nonumber \\
&\le\sum^{\el_0}_{k=0}\el_0!M_0\cdot\frac{R_0^{k+m_{j_0}-m_{i_0}}}{(k+m_{j_0}-m_{i_0})!}  \label{eq11}
\end{align}
We have
\[
m_{j_0}-m_{i_0}\ge m_{i_0+1}-m_{i_0}>N_1\ge N_0 \ \ \text{for} \ \ j_0\ge i_0+1.
\]
So we have $k+m_{j_0}-m_{i_0}>N_0$ for every $k=0,1,m\ld,\el_0$.

But we have $M_0\cdot\el_0!$ $\ga_v<1$ for every $v\ge N_0$ where $\ga_v=\dfrac{(2R_0)^v}{v!}$.\ This gives that
\[
M_0\cdot\el_0!\cdot\frac{(R_0)^{k+m_{j_0}-m_{i_0}}}{(k+m_{j_0}-m_{i_0})!}<1 \ \ \text{for every} \ \ k=0,1,\ld,\el_0.
\]
Thus
\[
M_0\el_0!\cdot\frac{R_0^{k+m_{j_0}-m_{i_0}}}{(k+m_{j_0}-m_{i_0})!}<
\frac{1}{2^{k+m_{j_0}-m_{i_0}}} \ \ \text{for every} \ \ k=0,1,\ld,\el_0.
\]
So we get
\begin{align}
&\sum^{\el_0}_{k=0}M_0\el_0!\frac{R_0^{k+m_{j_0}-m_{i_0}}}
{(k+m_{j_0}-m_{i_0})!}<\sum^{\el_0}_{k=0}\frac{1}{2^{k+m_{j_0}-m_{i_0}}} \nonumber \\
&\Rightarrow\sum^{\el_0}_{k=0}M_0\el_0!\frac{R_0^{k+m_{j_0}-m_{i_0}}}
{(k+m_{j_0}-m_{i_0})!}<\sum^{+\infty}_{v=m_{j_0}-m_{i_0}}\frac{1}{2^v} \nonumber \\
&=\frac{1}{2^{m_{j_0}-m_{i_0}}}\cdot\frac{1}{1-\dfrac{1}{2}}=\frac{1}
{2^{m_{j_0}-m_{i_0}-1}}.  \label{eq12}
\end{align}
By (\ref{eq11}) and (\ref{eq12}) we have:
\begin{eqnarray}
\|T_{m_{i_0},\la_0}(f_{j_0})\|_{\oD_{r_0}}<\frac{1}{2^{m_{j_0}-m_{i_0}-1}}. \label{eq13}
\end{eqnarray}
By (\ref{eq13}) we get:
\begin{eqnarray}
\|T_{m_{i_0},\la_{0}}(f_i)\|_{\oD_{R_0}}<\frac{1}{2^{m_j-m_{i_0}-1}} \ \ \text{for every} \ \ j=i_0+1,\ld,v_0.  \label{eq14}
\end{eqnarray}
By (\ref{eq14}) we obtain:
\begin{align*}
\sum^{v_0}_{i=i_0+1}\|T_{m_{i_0},\la_0}(f_i)\|_{\oD_{R_0}}&<
\sum^{v_0}_{j=i_0+1}\frac{1}{2^{m_j-m_{i_0}-1}}=2^{m_{i_0}+1}\sum^{v_0}_{j=i_0+1}
\frac{1}{2^{m_j}}\\
&<2^{m_{i_0}+1}\sum^{+\infty}_{v=m_{i_0+1}}\frac{1}{2^v}
=2^{m_{i_0}+1}\frac{1}{2^{m_{i_0+1}}}\cdots\frac{1}{1-\dfrac{1}{2}}\\
&=\frac{1}{2^{m_{i_0+1}-m_{i_0}-1}}.
\end{align*}
So we have:
\begin{eqnarray}
\sum^{v_0}_{i=i_0+1}\|T_{m_{i_0},\la_0}(f_i)\|_{\oD_{R_0}}<\frac{1}{2^{m_{i_0+1}
-m_{i_0}-2}}.  \label{eq15}
\end{eqnarray}
Now, by (\ref{eq7}) and (\ref{eq15}) we get:
\[
\|T_{m_{i_0},\la_0}(\vPi)-p\|_{\oD_{R_0}}<\|T_{m_{i_0},\de_{i_0}}(f_{i_0})-T_{m_{i_0},\la_0}(f_{i_0})\|_{\oD_{R_0}}+
\frac{1}{2^{m_{i_0+1}
-m_{i_0}-2}}
\]
and the proof of this lemma is complete now. \qb
\end{Proof}
\begin{lem}\label{lem2.7}
Let $(k_n)$ be a subsequence of natural numbers such that $\ssum^{+\infty}_{n=1}\dfrac{1}{k_n}=+\infty$.\ Then, for every positive number $M>0$ there exists a subsequence $(\mi_n)$ of $(k_n)$ such that
\begin{enumerate}
\item[(i)] $\mi_1>M$,
\item[(ii)] $\mi_{n+1}-\mi_n>M$ for every $n=1,2,\ld$ and
\item[(iii)] $\ssum^{+\infty}_{n=1}\dfrac{1}{m_n}=+\infty$.
\end{enumerate}
\end{lem}

After the previous lemmas we are ready now to prove Lemma \ref{lem2.3}.\vspace*{0.2cm}\\
\noindent
{\bf Proof of lemma \ref{lem2.3}}. We fix $n_0,j_0,s_0>2$, $\rho_0>2$ and we will prove that the set $\bigcup\limits^{+\infty}_{m=1}E(n_0,\rho_0,j_0,s_0,m)$ is dense in $(\ch(\C),\ct_u)$.

Let some polynomial $Q$, and $\e_1>0$. We consider the neighbourhood of $Q$,
\[
\cu(Q,K,\e_1):=\{h\in\ch(\C)\mid\|Q-h\|_K<\e_1\}
\]
for some fixed compact set $K\subseteq\C$. We denote $E:=\bigcup\limits^{+\infty}_{m=1}E(n_0,\rho_0,j_0,s_0,m)$ for simplicity.

It suffices to prove that $E\cap\cu(Q,K,\e_1)\neq\emptyset$.

We fix some positive number $R_0>1$ such that if we set $\oD_{R_0}:=\{z\in\C\mid |z|\le R_0\}$ we have $K\cup C_{n_0}\subset\oD_{R_0}$.

It suffices to find some polynomial $f$ and some natural number $m_0$ such that:

a) $\|f-Q\|_{\oD_{R_0}}<\e_0$ and

b) $f\in E(n_0,\rho_0,j_0,s_0,m_0)$, where $\e_0:=\min\Big\{\e_1,\dfrac{1}{s_0}\Big\}$.

Let $p_{j_0}(z)=\ssum^{\el_0}_{j=0}\bi_jz^j$, $\bi_j\in\C$, for $j=0,1,\ld,\el_0$ for some $\el_0\in\{0,1,2,\ld\}$.\ We set
\[
M_0:=\max\{|\bi_j|,\;j=0,1,\ld,\el_0\} \ \  \text{and}  \ \ M_1:=M_0\cdot\sum^{\el_0}_{j=0}R^j_0.
\]
Let some
\[
\de_0\in\bigg(0,\frac{1}{\rho_0}\log\bigg(1+\frac{\e_0}{4M_1}\bigg)\bigg)
\]
(for example $\de_0=\dfrac{1}{2\rho_0}\log\Big(1+\dfrac{\e_0}{4M_1}\Big)\Big)$.

Then we have:
\setcounter{equation}{0}
\begin{align}
&\de_0<\frac{1}{\rho_0}\log\bigg(1+\frac{\e_0}{4M_1}\bigg)\Rightarrow
\rho_0\de_0<\log\bigg(1+\frac{\e_0}{4M_1}\bigg)\Rightarrow \nonumber\\
&e^{\rho_0\de_0}<1+\frac{\e_0}{4M_1}\Leftrightarrow e^{\frac{\frac{\de_0}{1}}{\rho_0}}<1+\frac{\e_0}{4M_1}.  \label{eq1}
\end{align}

It is well known that
\[
\lim_{v\ra+\infty}\bigg(1+\frac{\rho_0\de_0}{v}\bigg)^v=e^{\rho_0\de_0}. \ \ \text{So}
\]
\begin{eqnarray}
\lim_{v\ra+\infty}\bigg(1+\frac{\rho_0\de_0}{k_v}\bigg)^{k_v}=e^{\rho_0\de_0}. \label{eq2}
\end{eqnarray}
By (\ref{eq1}) and (\ref{eq2}) we get that there exists some natural number $v_1\in\N$ such that:
\[
\bigg(1+\frac{\rho_0\de_0}{k_v}\bigg)^{k_v}<1+\frac{\e_0}{4M_1} \ \ \text{for every} \ \ v\in\N, \ \ v\ge v_1.
\]
So we have:
\begin{align}
&1+\frac{\rho_0\de_0}{k_v}<\sqrt[k_v]{1+\frac{\e_0}{4M_1}}\Rightarrow
\frac{\rho_0\de_0}{k_v}<\sqrt[k_v]{1+\frac{\e_0}{4M_1}}-1\nonumber \\
&\Rightarrow
\frac{\de_0}{k_v}<\frac{1}{\rho_0}\bigg(\sqrt[k_v]{1+\frac{\e_0}{4M_1}}-1\bigg) \ \ \text{for every} \ \  v\in\N, \ \ v\ge v_1. \label{eq3}
\end{align}
We consider the sequence
\[
\ga'_n:=\frac{(2\rho_0R_0)^v}{v!}\cdot\el_0!\cdot M_0, \ \ v=1,2,\ld\;.
\]
By ratio criterion for the convergence of sequences we conclude that $\ga'_n\ra0$ as $n\ra+\infty$, so there exists some natural number $v_2\in\N$ such that:
\begin{eqnarray}
\ga'_n<1 \ \ \text{for every} \ \ n\in\N, \ \ n\ge v_2.  \label{eq4}
\end{eqnarray}
We set
\[
v_3:=\max\bigg\{v_0,v_1,v_2,\el_0,\deg Q,3+\frac{1}{\log 2}\log\bigg(\frac{1}{\e_0}\bigg)\bigg\}+1
\]
where $\Big(1+\dfrac{\e_0}{2M_1}\Big)^{\frac{v_0}{v_0+\el_0}}>1+\dfrac{\e_0}{4M_1}$.

We apply \ref{lem2.7} and we get that there exists some subsequence $(\mi_n)$, $n=1,2,\ld$ of $(k_n)$ such that
\begin{enumerate}
\item[(i)] $\mi_1>v_3$
\item[(ii)] $\mi_{n+1}-\mi_n>v_3$ for every $n=1,2,\ld$
\item[(iii)] $\ssum^{+\infty}_{n=1}\dfrac{1}{\mi_n}=+\infty$.
\end{enumerate}
Now we have $\ssum^{+\infty}_{n=1}\dfrac{1}{\mi_n}=+\infty$.\ So $\ssum^{+\infty}_{n=1}\dfrac{\de_0}{\mi_n}=+\infty$.\ Thus, there exists the minimum natural number $N_0$ such that
\begin{eqnarray}
\sum^{N_0+1}_{n=1}\frac{\de_0}{\mi_n}>\rho_0-\frac{1}{\rho_0}.  \label{eq5}
\end{eqnarray}
By the above notations and terminology we define now a partition $\cp$ of the closed interval $\Big[\dfrac{1}{\rho_0},\rho_0\Big]$ as follows: We set
\[
a_1:=\frac{1}{\rho_0}, \ \ a_2:=\frac{1}{\rho_0}+\frac{\de_0}{\mi_1}, \ \ a_3:=a_2+\frac{\de_0}{\mi_2},\ld,a_{i+1}=a_i+\frac{\de_0}{\mi_i},
\]
for $i=1,2,\ld,N_0$, $a_{N_0+2}=\rho_0$ if $a_{N_0+1}<\rho_0$, or else  $a_{N_0+1}=\rho_0$ if the final point of the partition $\cp$ is $a_{N_0+1}$.

That is we have:
\[
\cp:=\bigg\{a_1=\frac{1}{\rho_0}<a_2<\cdots<a_{N_0}<a_{N_0+1}<a_{N_0+2}=\rho_0\bigg\}
\]
if $a_{N_0+1}<\rho_0$ or else
\[
\cp:=\bigg\{a_1=\frac{1}{\rho_0}<a_2<\cdots<a_{N_0+1}=\rho_0\bigg\}, \ \ \text{if} \ \  a_{N_0+1}=\rho_0
\]
where $a_1=\dfrac{1}{\rho_0}$, $a_{i+1}=a_i+\dfrac{\de_0}{\mi_i}$ for $i=1,2,\ld,N_0$.

After the above preparation, we are ready now to define the function $f\in\ch(\C)$ that satisfies the properties (a), (b) above.

We consider the solutions $(\mi_i,a_i,p_{j_0})$ for $i=1,2,\ld,N_0$ if $a_{N_0+1}<\rho_0$ and the solutions $(\mi_i,a_i,p_{j_0})$ for $i=1,2,\ld,N_0$ and the solution $(\mi_{N_0+1},\rho_0,p_{j_0})$ also if $a_{N_0+1}=\rho_0$.

We denote the solution $(\mi_i,a_i,p_{j_0})$ as $f_i$, for $i=1,2,\ld,N_0$ if $a_{N_0+1}<\rho_0$ and as $f_{n_0+1}$ the solution $(\mi_{N_0+1},\rho_0,p_{j_0})$ if $a_{N_0+1}=\rho_0$.

Now we set: $f:=Q+\ssum^{N_0}_{i=1}f_i$ if $a_{N_0+1}<\rho_0$ and $f:=Q+\ssum^{N_0+1}_{i=1}f_i$ if $a_{N_0+1}=\rho_0$.

Obviously, $f$ is a polynomial, so $f\in\ch(\C)$.

We prove now that $f$ satisfies the properties (a) and (b).

Firstly, we prove that $f$ satisfies property (a).

We suppose that $a_{N_0+1}<\rho_0$.\ In this case we have:
\begin{eqnarray}
\|f-Q\|_{\oD_{R_0}}=\bigg\|\bigg(Q+\sum^{N_0}_{i=1}f_i\bigg)-Q\bigg\|_{\oD_{R_0}}
=\bigg\|\sum^{N_0}_{i=1}f_i\bigg\|_{\oD_{R_0}}\le\sum^{N_0}_{i=1}\|f_i\|_{\oD_{R_0}}. \label{eq6}
\end{eqnarray}
The polynomial $f_i$, for $i=1,2,\ld,N_0$ is the $(\mi_i,a_i,p_{j_0})$ solution, that is
\[
f_i(z)=\sum^{\el_0}_{k=0}\frac{k!}{(k+\mi_i)!}\frac{\bi_k}{a^{k+m_i}_i}
z^{k+\mi_i} \ \ \text{for} \ \ i=1,2,\ld,N_0.
\]
So, for $i=1,2,\ld,N_0$, $z\in\oD_{R_0}$ we have:
\begin{align}
|f_i(z)|&=\bigg|\sum^{\el_0}_{k=0}\frac{k!}{(k+\mi_i)!}\frac{\bi_k}{a^{k+\mi_i}_i}z^{k+\mi_i}\bigg| \le\sum^{\el_0}_{k=0}\frac{k!}{(k+\mi_i)!}\frac{|\bi_k|}{a^{k+\mi_i}_i}R^{k+\mi_i}_0 \nonumber \\
&\le\el_0!\cdot M_0\cdot\sum^{\el_0}_{k=0}\bigg(\frac{R_0}{a_i}\bigg)^{k+\mi_i}
\frac{1}{(k+\mi_i)!}\le\el_0!\cdot M_0\cdot\sum^{\el_0}_{k=0}\frac{(\rho_0R_0)^{k+\mi_i}}{(k+\mi_i)!}.  \label{eq7}
\end{align}
For every $k=0,1,\ld,\el_0$ we have $k+\mi_i\ge\mi_1>v_3>v_2$.\ So by (\ref{eq4}) we get
\begin{align}
&\frac{(2\rho_0R_0)^{k+\mi_i}}{(k+\mi_i)!}\el_0!\cdot M_0<1 \ \ \text{for every} \ \ k=0,1,\ld,\el_0\;\Rightarrow \nonumber \\
&\el_0!\cdot M_0\frac{(\rho_0R_0)^{k+\mi_i}}{(k+\mi_i)!}<\frac{1}{2^{k+\mi_i}} \ \ \text{for every} \ \ k=0,1,\ld,\el_0\;\Rightarrow \nonumber\\
&\sum^{\el_0}_{k=0}\el_0!\cdot M_0\frac{(R_0\rho_0)^{k+\mi_i}}{(k+\mi_i)!}<\sum^{\el_0}_{k=0}\frac{1}{2^{k+\mi_i}}<
\sum^{+\infty}_{v=\mi_i}\frac{1}{2^v}=\frac{1}{2^{\mi_i-1}}.  \label{eq8}
\end{align}
By (\ref{eq7}) and (\ref{eq8}) we get
\begin{eqnarray}
\|f_i\|_{\oD_{R_0}}<\frac{1}{2^{\mi_i-1}} \ \ \text{for every} \ \ i=1,2,\ld,N_0.  \label{eq9}
\end{eqnarray}
By (\ref{eq9}) we have:
\begin{eqnarray}
\sum^{N_0}_{i=1}\|f_i\|_{\oD_{R_0}}<\sum^{N_0}_{i=1}\frac{1}{2^{\mi_i-1}}<
\sum^{+\infty}_{v=\mi_1-1}\frac{1}{2^v}=\frac{1}{2^{\mi_1-2}}.  \label{eq10}
\end{eqnarray}
By (\ref{eq6}) and (\ref{eq10}) we obtain:
\[
\|f-Q\|_{\oD_{R_0}}<\frac{1}{2^{\mi_1-2}}<\e_0
\]
because $\mi_1>v_3$ and so $\mi_1>2+\dfrac{1}{\log2}\cdot\log\Big(\dfrac{1}{\e_0}\Big)$ by hypothesis (i) for $\mi_1$.

So we have proved that $\|f-Q\|_{\oD_{R_0}}<\e_0$ that is property (a) for the function $f$.\ It remains to prove b).

We prove b) for $m_0=\mi_{N_0}$.\ That is we prove now that $f\in E(n_0,\rho_0,j_0,s_0,\mi_{N_0})$.

We fix some $\la_0\in\De_{\rho_0}=\Big[\dfrac{1}{\rho_0},\rho_0\Big]$.\ We consider the unique $i_0\in\{1,2,\ld,N_0+1\}$ such that $\la_0\in[a_{i_0},a_{i_0+1})$ if $1\le i_0\le N_0$ or we set $i_0=N_0+1$ if $\la_0\in[a_{N_0+1},\rho_0]$.\ We show now that:
\begin{eqnarray}
\|T_{\mi_{i_0},\la_0}(f)-p_{j_0}\|_{\oD_{R_0}}<\frac{1}{s_0}.  \label{eq11}
\end{eqnarray}
We remark that all the suppositions of lemma are satisfied, so we get:
\begin{eqnarray}
\|T_{\mi_{i_0},\la_0}(f)-p_{j_0}\|_{\oD_{R_0}}<\|T_{\mi_{i_0},a_{i_0}}(f_{i_0})-
T_{\mi_{i_0},\la_0}(f_{i_0})\|_{\oD_{R_0}}+\frac{1}{2^{\mi_{i_0+1}-(\mi_{i_0}+2)}}. \label{eq12}
\end{eqnarray}
We show that:
\[
\la_0\in\bigg[a_{i_0},a_{i_0}\cdot\sqrt[\mi_{i_0}+\el_0]{1+\e_0/2M_1}\bigg).
\]
We have that: $\la_0\in[a_{i_0},a_{i_0+1})$ by hypothesis.\ It suffices to have:
\begin{align}
a_{i_0+1}<a_{i_0}\sqrt[\mi_{i_0}+\el_0]{1+\frac{\e_0}{2M_1}}&\Leftrightarrow a_{i_0+1}-a_{i_0}<a_{i_0}\bigg(\sqrt[\mi_{i_0}+\el_0]{1+\frac{\e_0}{2M_1}}-1\bigg)
\nonumber \\
&\Leftrightarrow\frac{\de_0}{\mi_{i_0}}<a_{i_0}\cdot\bigg(\sqrt[\mi_{i_0+\el_0}]
{1+\frac{\e_0}{2M_1}}-1\bigg).  \label{eq13}
\end{align}
It suffices to have:
\begin{eqnarray}
\frac{\de_0}{\mi_{i_0}}<\frac{1}{\rho_0}\bigg(\sqrt[\mi_{i_0}+\el_0]{1+\frac{\e_0}{2M_1}}-1\bigg)
\label{eq14}
\end{eqnarray}
because $\dfrac{1}{\rho_0}\le a_{i_0}$.

By our hypothesis we have:
\begin{eqnarray}
\frac{\de_0}{\mi_{i_0}}<\frac{1}{\rho_0}\cdot\bigg(\sqrt[\mi_{i_0}]
{1+\frac{\e_0}{4M_1}}-1\bigg).  \label{eq15}
\end{eqnarray}
So, in order (\ref{eq14}) holds it suffices by (\ref{eq15}) to hold:
\begin{align}
&\frac{1}{\rho_0}\bigg(\sqrt[\mi_{i_0}]{1+\frac{\e_0}{4M_1}}-1\bigg)<\frac{1}{\rho_0}
\cdot\bigg(\sqrt[\mi_{i_0}+\el_0]{1+\frac{\e_0}{2M_1}}-1\bigg)\nonumber \\
&\Leftrightarrow1
+\frac{\e_0}{4M_1}<\bigg(1+\frac{\e_0}{2M_1}\bigg)^{\frac{\mi_{i_0}}{\mi_{i_0}+\el_0}}.
\label{eq16}
\end{align}

We have
\begin{eqnarray}
\lim_{x\ra1^-}\bigg(1+\frac{\e_0}{2M_1}\bigg)^x=1+\frac{\e_0}{2M_1}>1+\frac{\e_0}{4M_1}.
\label{eq17}
\end{eqnarray}
So, by (\ref{eq17}) inequality (\ref{eq16}) holds, because $\mi_{i_0}\ge\mi_1>v_3>v_0$ and the definition of $v_0$.

Thus, we have
\[
\la_0\in[a_{i_0},a_{i_0}\cdot\sqrt[\mi_{i_0}+\el_0]{1+\frac{\e_0}{2M_1}}\bigg)
\]
and by Lemma \ref{lem2.5} we get:
\begin{eqnarray}
\|T_{\mi_{i_0},a_{i_0}}(f_{i_0})-T_{\mi_{i_0},\la_0}(f_{i_0})\|_{\oD_{R_0}}<
\frac{\e_0}{2}.  \label{eq18}
\end{eqnarray}
By the hypothesis $\mi_1>v_3$ we get
\begin{eqnarray}
\frac{1}{2^{\mi_{i_0+1}-(\mi_{i_0}+2)}}<\frac{\e_0}{2}.  \label{eq19}
\end{eqnarray}
Finally, by (\ref{eq12}), (\ref{eq18}) and (\ref{eq19}) we get:
\[
\|T_{\mi_{i_0},\la_0}(f)-p_{j_0}\|_{\oD_{R_0}}<\e_0\le\frac{1}{s_0}
\]
and we have proved inequality (\ref{eq11}).\ The other case where $i_0=N_0+1$ is proved with nearly almost the same way.

The proof of this lemma is complete now. \qb \vspace*{0.2cm}

Now by Lemmas \ref{lem2.1}, \ref{lem2.2}, \ref{lem2.3}, the fact that the space $(\ch(\C),\ct_u)$ is a complete metric space and Baire's Category Theorem the proof of Proposition \ref{prop1.1} is complete.
\section{The general case: proof of Theorem \ref{thm1.4}}\label{sec3}
The notation introduced in the previous section will be used in the present section as well. For instance, the sequence $(p_j)$ denotes an enumeration of all the polynomials with coefficients in $\mathbb{Q}+i\mathbb{Q}$.
We introduce some extra terminology and a non-trivial result from the book \cite{KN}.

\noindent
{\bf Uniform distribution Modulo 1} \smallskip

For a real number $x$, let $[x]$ denote the integer part of $x$, that is, the greatest integer $\le x$; let $\{x\}=x-[x]$ be the fractional part of $x$, or the residue of $x$ modulo 1. We note that the fractional part of any real number is contained in the unit interval $I=[0,1)$.

Let $\oo=(x_n)$ $n=1,2,\ld$, be a given sequence of real numbers. For a positive integer $N$ and a subset $E$ of $I$, let the counting function $A(E;N;\oo)$ be defined as the number of terms $x_n$, $1\le n\le N$, for which $\{x_n\}\in E$.

If there is no risk of confusion, we shall often write $A(E;N)$ instead of $A(E;N;\oo)$. Here is our basic definition.
\begin{Def}\label{Def3.3}
The sequence $\oo=(x_n)$, $n=1,2,\ld$, of real numbers is said to be uniformly distributed modulo 1 (abbriviated u.d. $\mod 1$) if for every pair $a,b$ of real numbers with $0\le a<b\le1$ we have
\[
\lim_{N\ra+\infty}\frac{A([a,b);N;\oo)}{N}=b-a.
\]
Thus, in simple terms, the sequence $(x_n)$ is u.d. $\mod1$ if every half-open subinterval of $I$ eventually gets its ``proper share'' of fractional parts.
\end{Def}

Let $(\cu_n(x))$, $n=1,2,\ld$, for every $x$ lying in some given bounded or unbounded interval $J$, be a sequence of real numbers. The sequence $(\cu_n(x))$ is said to be u.d. $\mod1$ for almost all $x$, if for every $x\in J$, apart from a set which has Lebesgue measure 0, the sequence $(\cu_n(x))$ is u.d. $\mod1$. The proof of Theorem \ref{thm1.4} uses in a crucial way the following theorem, due to Weyl \cite{KN}.
\begin{thm}(Weyl)\label{thm3.4} 
Let $(a_n)$, $n=1,2,\ld$, be a given sequence of distinct integers. Then the sequence $(a_nx)$, $n=1,2,\ld$, is u.d. $\mod1$ for almost all real numbers $x$.
\end{thm}
We fix some further notation.
Let $\la>0$ and some $f\in\ch C(\{T_{k_v,\la}\})$.\ We can fix some sequence $(\mi^j_n)$ from distinct integers, $n=1,2,\ld$, for every $j=1,2,\ld$ such that
\[
T_{\mi^j_n,\la}(f)\ra p_j.
\]

To avoid confusion, we stress that in the sequence $(\mi^j_n)$, $j$ is just an index and \textit{not a power}. This will be the case whenever the sequence $(\mi^j_n)$ occurs throughout this section.\\
By Theorem \ref{thm3.4} we have that the sequence $(\thi\mi^j_n)$ is u.d $\mod1$ for almost all $\thi\in[0,1)$ for every $j=1,2,\ld\;.$

We set
\[
\ca_{(\mi^j_n)}:=\{\thi\in[0,1)\mid\text{the sequence}\;(\thi\mi^j_n),\ \
n=1,2,\ld \ \ \text{is not u.d}  \mod1\}
\]
for every $j=1,2,\ld\;.$

Then $\la_1(\ca_{(\mi^j_n)})=0$ for every $j=1,2,\ld$.\ We set $\ca':=\bigcup\limits^{+\infty}_{j=1}\ca_{(\mi^j_n)}$.\ Then $\la_1(\ca')=0$.\ We denote $\ca:=[0,1)\sm\ca'$. Then
\[
\ca:=\{\thi\in[0,1)\mid\text{the sequence}\;(\thi\mi^j_n),\;
n=1,2,\ld \; \text{is u.d}  \mod1 \; \text{for every} \; j=1,2,\ld\}
\]
and $\la_1(\ca)=1$.

We consider the set:
\[
\cb_{\la}:=\{\thi\in[0,1)\mid f\in\ch C(\{T_{k_v,\la e^{2\pi i\thi}}\})\}.
\]
We will prove the following proposition.
%
\begin{prop}\label{prop3.6}
By the above notations and terminology we have: $\ca\subseteq\cb_\la$.
\end{prop}

Before we prove the above Proposition \ref{prop3.6}, we will prove some helpful lemma.

We fix some $\thi_0\in\ca$, $j_0\in\N$, $n_0>2$, $n_0\in\N$, $\la_0>0$ and $\e_0\in(0,1)$.\ We consider the polynomial $g(z):=p_{j_0}(e^{2\pi i\thi_0}z)$.
\begin{lem}\label{lem3.7}
There exists some $v_0\in\N$ such that:
\[
\|T_{k_{v_0},\la_0e^{2\pi i\thi_0}}(f)-g\|_{\oD_{n_0}}<\e_0.
\]
\end{lem}
\begin{Proof}
We set $M_0:=\|p_{j_0}\|_{\oD_{n_0}}$.

We consider the trinomial $A(x):=x^2+(M_0+1)x-\e_0$.\ Because $\e_0\in(0,1)$, this trinomial has two roots $\rho_1$ and $\rho_2$ such that $\rho_1<0<\rho_2$.\ We choose some positive number $\e_1\in(0,\rho_2)$, for example take $\e_1:=\dfrac{\rho_2}{2}$.\ Then we have:
\setcounter{equation}{0}
\begin{eqnarray}
\e^2_1+(M_0+1)\e_1<\e_0\Leftrightarrow A(\e_1)<0.  \label{eq1}
\end{eqnarray}
We solve now the inequality $|e^{2\pi i\thi_0v}-1|<\e_1$ in the set of natural numbers.

Let
\[
L:=\{v\in\N\mid |e^{2\pi i\thi_0v}-1|<\e_1\}.
\]
We compute easily that
\begin{eqnarray}
|e^{2\pi i\thi_0v}-1|=2\cdot|\sin(\pi\thi_0v)|. \label{eq2}
\end{eqnarray}
So, by (\ref{eq2}) we get
\[
L=\bigg\{v\in\N\mid|\sin(\pi\thi_0v)|<\frac{\e_1}{2}\bigg\}.
\]
Let $\f_0$ be the unique arc such that
\begin{eqnarray}
\f_0\in\bigg(0,\frac{\pi}{2}\bigg) \ \ \text{and} \ \ \sin\f_0=\frac{\e_1}{2}.  \label{eq3}
\end{eqnarray}
Using elementary trigonometry we can see that
\begin{eqnarray}
L':=\bigg\{v\in\N\mid\{\thi_0v\}\in\bigg(0,\frac{\f_0}{\pi}\bigg)\cup\bigg(
1-\frac{\f_0}{\pi},1\bigg)\bigg\}\subseteqq L.  \label{eq4}
\end{eqnarray}
We have:
\[
T_{\mi^{j_0}_n,\la_0}(f)\ra p_{j_0} \ \ \text{as} \ \  n\ra+\infty,
\]
so there exists $v_1\in\N$ such that for every $n\in\N$, $n\ge v_1$ holds:
\begin{eqnarray}
\|T_{\mi^{j_0}_n,\la_0}(f)-p_{j_0}\|_{\oD_{n_0}}<\e_1.  \label{eq5}
\end{eqnarray}
We can see easily that for every $v\in\N$ we obtain that:
\begin{align}
\|T_{k_v,\la_0}e^{2\pi i\thi_0}-g\|_{\oD_{n_0}}\le&\,|e^{2\pi i\thi_0k_v}-1|\cdot(\|T_{k_v,\la_0}(f)-p_{j_0}\|_{\oD_{n_0}}+\|p_{j_0}\|_{\oD_{n_0}}) \nonumber \\
&+\|T_{k_v,\la_0}(f)-p_{j_0}\|_{\oD_{n_0}}.  \label{eq6}
\end{align}
By (\ref{eq5}) and (\ref{eq6}) we get that for every $n\in\N$, $n\ge v_1$ holds:
\begin{eqnarray}
\|T_{\mi^{j_0}_n,\la_0}e^{2\pi i\thi_0}-g\|_{\oD_{n_0}}<
|e^{2\pi i\thi_0\mi^{j_0}_n}-1|\cdot(\e_1+M_0)+\e_1.  \label{eq7}
\end{eqnarray}
We suppose now that there exists not any natural number $n>v_1$ such that
\[
\{\thi_0\mi^{j_0}_n\}\in\bigg(0,\frac{\f_0}{\pi}\bigg)\cup\bigg(1-\frac{\f_0}{\pi},1\bigg).
\]
So, we have $\{\thi_0\mi^{j_0}_n\}\notin\Big(0,\dfrac{\f_0}{\pi}\Big)$, for every $n>v_1$.

Then, by definition we have:
\[
A\bigg(\bigg[\frac{\f_0}{2\pi},\frac{\f_0}{\pi}\bigg),N,(\thi_0\mi^{j_0}_n)\bigg)\le
\mi^{j_0}_{v_1} \ \ \text{for every} \ \ n\in\N,
\]
thus:
\[
\lim_{N\ra+\infty}\frac{A\Big(\Big[\dfrac{\f_0}{2\pi},\dfrac{\f_0}{\pi}\Big),N,(\thi_0\mi^{j_0}_n)\Big)}
{N}=0
\]
that is not true because the sequence $(\thi_0\mi^{j_0}_n)$ is u.d. $\mod1$ by its definition, and so
\[
\lim_{N\ra+\infty}\frac{A\Big(\big[\dfrac{\f_0}{2\pi},\dfrac{\f_0}{\pi}\Big),N,(\thi_0\mi^{j_0}_n)\Big)}{N}=
\frac{\f_0}{2\pi}.
\]
So, there exists some natural number $n_1>v_1$ such that
\begin{eqnarray}
\{\thi_0\mi^{j_0}_{n_1}\}\in\bigg(0,\frac{\f_0}{\pi}\bigg)\cup\bigg(1-\frac{\f_0}{\pi},1\bigg).
\label{eq8}
\end{eqnarray}
By (\ref{eq4}) and (\ref{eq8}) we have
\begin{eqnarray}
\mi^{j_0}_{n_1}\in L,  \label{eq9}
\end{eqnarray}
thus by the definition of $L$ we get:
\begin{eqnarray}
|e^{2\pi i\thi_0\mi^{j_0}_{n_1}}-1|<\e_1.  \label{eq10}
\end{eqnarray}
Finally, by (\ref{eq1}), (\ref{eq7}) and (\ref{eq10}) we have:
\[
\|T_{\mi^{j_0}_{n_1},\la_0}e^{2\pi i\thi_0}-g\|_{\oD_{n_0}}<\e_0,
\]
that is for the unique $v_0\in\N$ such that $k_{v_0}=\mi^{j_0}_{n_1}$ we have the conclusion of this\linebreak lemma. \qb
\end{Proof}

After the above Lemma \ref{lem3.7} we are ready now to prove Proposition \ref{prop3.6}.\vspace*{0.2cm} \\
\noindent
{\bf Proof of Proposition \ref{prop3.6}.} Let fixed $\thi_0\in\ca$. Because the set $\vPsi:=\{p_1,p_2,\ld\}$ is dense, we conclude easily that the set
\[
\vPsi_{\thi_0}:=\{p_1(e^{2\pi i\thi_0}z),\;p_2(e^{2\pi i\thi_0}z),\ld,\}=
\{p_j(e^{2\pi i\thi_0}z),\;j=1,2,\ld,\}
\]
is dense.\ By this fact and Lemma \ref{lem3.7} we conclude that the set $\{T_{k_v,\la_0}e^{2\pi i\thi_0}(f)$, $v=1,2,\ld\}$ is dense, that is $f\in\ch C\{T_{k_v,\la_0}e^{2\pi i\thi_0}\})$ and so $\thi_0\in\cb_{\la_0}$.\ The proof of Proposition \ref{prop3.6} is complete now. \qb \vspace*{0.2cm}

We fix now some
\[
f\in\ch=\bigcap_{\la\in(0,+\infty)}\ch C(\{T_{k_v,\la}\})
\]

We consider the set
\[
\widetilde{J}(f):=\{z\in\C\sm\{0\}\mid f\in\ch C(\{T_{k_n,z}\})\}.
\]
Of course, by the definition of the set $\widetilde{J}(f)$ we have: $(0,+\infty)\subseteq\widetilde{J}(f)$.\ Now, by the respective definitions we can see easily that:
\[
\widetilde{J}(f)_\la=\cb_\la \ \ \text{for every} \ \ \la>0.
\]
Applying Proposition \ref{prop3.6} we conclude Theorem \ref{thm1.4} for $J=\widetilde{J}(f)$ and this completes our work.

{\bf Acknowledgment} I would like to thank G. Costakis for all the help he offered me concerning this project.

Department of Mathematics and Applied Mathematics, University of Crete, Panepistimiopolis Voutes, 700-13, Heraklion, Crete, Greece.\\
email:tsirivas@uoc.gr

\end{document}